\documentclass[leqno,11pt]{article}

\usepackage{amssymb,amsmath,amscd,amsfonts,a4,psfrag,graphicx,latexsym,epsfig,color}

\newtheorem{theo}{Theorem}[section]

\newtheorem{lem}[theo]{Lemma}

\newtheorem{prop}[theo]{Proposition}

\newtheorem{rem}[theo]{Remark}
\newtheorem{coro}[theo]{Corollary}

\newtheorem{exam}[theo]{Example}

\newcommand{\ugot}{\ensuremath{\mathfrak{u}}}
\newcommand{\kgot}{\ensuremath{\mathfrak{k}}}

\newcommand{\tgot}{\ensuremath{\mathfrak{t}}}

\newcommand{\pgot}{\ensuremath{\mathfrak{p}}}

\newcommand{\glgot}{\ensuremath{\mathfrak{gl}}}

\newcommand{\spgot}{\ensuremath{\mathfrak{sp}}}
\newcommand{\Rgot}{\ensuremath{\mathfrak{R}}}

\newcommand{\Ccal}{\ensuremath{\mathcal{C}}}

\newcommand{\Ocal}{\ensuremath{\mathcal{O}}}
\newcommand{\Pcal}{\ensuremath{\mathcal{P}}}

\newcommand{\Cbb}{\ensuremath{\mathbb{C}}}
\newcommand{\Rbb}{\ensuremath{\mathbb{R}}}

\newcommand{\horn}{\ensuremath{\hbox{\rm Horn}}}
\newcommand{\SHorn}{\ensuremath{{\rm Horn}_{\rm sp}}}

\newcommand{\trace}{\ensuremath{{\rm Tr}}}

\begin{document}

\title{The Horn cone associated with symplectic eigenvalues}

\author{Paul-Emile Paradan\footnote{IMAG, Univ Montpellier, CNRS, email : paul-emile.paradan@umontpellier.fr}}

\maketitle

\date{}


\begin{abstract}
In this note, we show that the Horn cone associated with symplectic eigenvalues admits the same inequalities as the classical Horn cone, except that the equality 
corresponding to $\trace(C)=\trace(A)+\trace(B)$ is replaced by the inequality corresponding to $\trace(C)\geq \trace(A)+\trace(B)$. 
\end{abstract}




\section{Introduction}

We consider $\Rbb^{2n}$ equipped with its canonical symplectic structure $\Omega_n=\sum_{k=1}^n d x_{k}\wedge dx_{k+n}$. Recall that a  family $(e_k)_{1\leq k\leq 2n}$ is a 
symplectic basis of $\Rbb^{2n}$, if $\Omega_n(e_{k}, e_{\ell})=0$ if $| k-\ell |\neq n$ and $\Omega_n(e_{k}, e_{k+n})=1$, $\forall k$.

Williamson's theorem \cite{Williamson} says that any positive definite quadratic form $q:\Rbb^{2n}\to \Rbb$ can be written $q(v)=\sum_{k=1}\lambda_k (v_{k}^2+v_{k+n}^2)$ 
where the $(v_j)$ are the coordinates of the vector $v\in\Rbb^{2n}$ relatively to a symplectic basis. The positive numbers $\lambda_k$, that one choose so that 
$$
\lambda(q):=(\lambda_1\geq \cdots\geq \lambda_n),
$$
will be referred to as the {\em symplectic eigenvalues} of the quadratic form $q$. They correspond to the frequencies of the normal modes of oscillation for the linear Hamiltonian system generated by $q$.

The object of study of this note concerns the symplectic Horn cone, denoted $\SHorn(n)$, that is defined as the set of triplets 
$(\lambda(q_1),\lambda(q_2),\lambda(q_1+q_2))$ where 
$q_1,q_2$ are positive definite quadratic forms on $\Rbb^{2n}$.

\begin{exam}In dimension $2$, the symplectic eigenvalue $\lambda(q)$ of a positive definite quadratic form $q(x_1,x_2)= a x_1^2 + b x_2^2+cx_1x_2$ is equal to $\tfrac{1}{2}\sqrt{4ab-c^2}$. 
It is straightforward to show that $\SHorn(1)$ is equal to the set of triplets $(x,y,z)$ of positive numbers satisfying $x+y\leq z$.
\end{exam}

 Our main Theorem states that $\SHorn(n)$ is a convex polyhedral set. Before detailing it, let us recall
some related results.

In \cite{Weinstein01}, A. Weinstein showed that for non-increasing $n$-tuples of positive real numbers $a$ and $b$, the set 
$\Delta_{\rm sp}(a,b):=\{\lambda(q_1+q_2)\, \vert\, \lambda(q_1)=a, \lambda(q_2)=b\}$ is  closed, convex and locally polyhedral. 

Recently, several authors have 
realized that some inequalities obtained long ago in the context of eigenvalues of Hermitian matrices still apply to symplectic eigenvalues:
\begin{itemize}
\item T. Hiroshima proved in \cite{Hiroshima} an analogue of Ky Fan inequalities: $\sum_{j=1}^k \lambda_j(q_1+q_2)\geq \sum_{j=1}^k \lambda_j(q_1)+\sum_{j=1}^k \lambda_j(q_2)$
\item In \cite{Jain-Mishra},  T. Jain and H. Mishra obtained an analogue of Lidskii inequalities: $\sum_{j=1}^k \lambda_{i_j}(q_1+q_2)\geq \sum_{j=1}^k \lambda_{i_j}(q_1)+\sum_{j=1}^k \lambda_j(q_2)$ for any subset 
$\{i_1<i_2<\cdots <i_k\}$.
\item In \cite{Bhatia-Jain}, R. Bhatia and T. Jain  obtained an analogue of the Weyl inequalities: $\lambda_{i+j-1}(q_1+q_2)\geq \lambda_{i}(q_1)+\lambda_{j}(q_2)$.
\end{itemize}

As the previous results suggest, we know explain the strong relationship between $\SHorn(n)$ with the classical Horn cone. If $A$ is a Hermitian $n\times n$ matrix, 
we denote by $\mathrm{s}(A)=(\mathrm{s}_1(A)\geq \cdots \geq\mathrm{s}_n(A))$ its spectrum. The Horn cone $\horn(n)$ is defined as the set of triplets 
$(\mathrm{s}(A),\mathrm{s}(B),\mathrm{s}(A+B))$ where $A,B$ are Hermitian $n\times n$ matrices.

Denote the set of cardinality $r$ subsets $I=\{i_1<i_2<\cdots<i_r\}$ of $[n]:=\{1,\ldots,n\}$ by $\Pcal^n_r$. To each $I\in \Pcal^n_r$
we associate: 
\begin{itemize}
\item a weakly decreasing sequence of non-negative integers 
$\lambda(I)=(\lambda_1\geq\cdots\geq \lambda_r)$ where $\lambda_a= n-r+a-i_a$ for $a\in [r]$.
\item the irreducible representation $V_{\lambda(I)}$ of $GL_r(\Cbb)$ with highest weight $\lambda(I)$.
\end{itemize}

If $x=(x_1,\ldots,x_n)\in\Rbb^n$ and $I\subset [n]$, we define $|\,x\,|_I=\sum_{i\in I}x_i$ and $|\,x\,|=\sum_{i=1}^n x_i$. Let us denote by $\Rbb_{+}^n$ the set of non-increasing $n$-tuples of real numbers.

\medskip

A. Klyachko \cite{Klyachko} has shown that an element $(x,y,z)\in (\Rbb_{+}^n)^3$ belongs to the cone $\horn(n)$ if and only if it satisfies $|x|+|y|=|z|$ and 
\begin{equation*}
|x|_I+|y|_J\leq |z|_K \tag{$(\star)_{I,J,K}$}
\end{equation*}
for any $r<n$, for any $I,J,K\in \Pcal^n_r$ such that the Littlewood-Richardson coefficient 
$$
c_{IJ}^{K}:= \dim [V_{\lambda(I)}\otimes V_{\lambda(J)}\otimes V_{\lambda(K)}^*]^{GL_r(\Cbb)}
$$
is non-zero. P. Belkale \cite{Belkale01} showed that the inequalities $(\star)_{I,J,K}$ associated to the condition $c_{IJ}^{K}=1$ are sufficient. Finally 
A. Knutson, T. Tao, and C. Woodward \cite{Knutson-Tao-Woodward-04}  have proved that this smaller list is actually minimal.  
We refer the reader to survey articles \cite{Fulton-00,Brion-Bourbaki} for details.

\newpage

The main result of this note is the following Theorem. Let us denote by $\Rbb_{++}^n$ the set of non-increasing $n$-tuples of positive real numbers.

\begin{theo}\label{theo:main}
An element $(x,y,z)\in(\Rbb^n_{++})^3$ belongs to $\SHorn(n)$ if and only if it satisfies
\begin{enumerate}
\item $|x|+|y|\leq |z|$,
\item $(\star)_{I,J,K}$ for all $(I,J,K)$ of cardinality $r<n$ such that  $c_{IJ}^{K}=1$.
\end{enumerate}
\end{theo}

\begin{coro}
Let $a,b\in \Rbb_{++}^n$. An element $z\in \Rbb_{++}^n$ belongs to $\Delta_{\rm sp}(a,b)$ if and only if it satisfies 
$|a|+|b|\leq |z|$ and $|a|_I+|b|_J\leq |z|_K$ for all $(I,J,K)$ of cardinality $r<n$ such that  $c_{IJ}^{K}=1$.
\end{coro}

\section{The causal cone of the symplectic Lie algebra}\label{sec:holomorphic}

The $2n\times 2n$ matrix $J_{n}= \begin{pmatrix}
0& -I_n\\
I_n& 0
\end{pmatrix}$ defines a complex structure on $\Rbb^{2n}$ that is compatible with the symplectic structure $\Omega_n$. The symplectic group 
$Sp(\Rbb^{2n})$ is defined by the relation ${}^t g J_n g=J_n$. A matrix $X$ belongs to the Lie algebra $\spgot(\Rbb^{2n})$ of 
$Sp(\Rbb^{2n})$ if and only the matrix $J_{n}X$ is symmetric. Moreover, $J_n X$ is positive if and only if 
$\Omega_n(Xv,v)\geq 0$, $\forall v\in \Rbb^{2n}$.

We call an invariant convex cone $C$ in $\spgot(\Rbb^{2n})$ a causal cone if $C$ is nontrivial, closed, and satisfies $C \cap -C=\{0\}$. A classical result \cite{Vinberg80,Paneitz81,Paneitz83} asserts that 
there are exactly two causal cones in $\spgot(\Rbb^{2n})$ : one, denoted $\mathbf{C}(n)$, containing $-J_n$ and its opposite $-\mathbf{C}(n)$. The causal cone $\mathbf{C}(n)$ is determined by 
the following equivalent conditions : for $X\in \spgot(\Rbb^{2n})$, we have 
$$
X\in \mathbf{C}(n)\hspace{3mm} \Longleftrightarrow\hspace{3mm} J_n X \ \text{is positive}\hspace{3mm}
\Longleftrightarrow\hspace{3mm} \trace(XgJ_ng^{-1})\geq 0,\, \forall g\in Sp(\Rbb^{2n}).
$$

Now we explain how is parameterized the interior $\mathbf{C}(n)^0$ of $\mathbf{C}(n)$. From the definition above, we see first that $X\in \mathbf{C}(n)^0$ if and only if $J_nX$ is positive definite.

The Lie algebra of the maximal compact subgroup $K=Sp(2n,\Rbb)\cap O(2n)$ is 
$$
\kgot:=
\left\{
\begin{pmatrix}A& B\\ -B& A  \end{pmatrix},\ {}^t A= -A, \ {}^t B= B
\right\}.
$$
If $\mu:=(\mu_1,\cdots,\mu_n)$, we write $\Delta(\mu)={\rm Diag}(\mu_1,\cdots,\mu_n)$ and 
$X(\mu)=\begin{pmatrix}0& \Delta(\mu)\\ -\Delta(\mu)& 0  \end{pmatrix}$. We work with the Cartan subalgebra $\tgot:=\{X(\mu),\mu\in\Rbb^n\}$ of $\kgot$ and the corresponding maximal torus 
$T\subset K$. The set of roots $\Rgot$ relatively to the action of $T$ on $\spgot(\Rbb^{2n})\otimes\Cbb$ are composed by the compact ones $\Rgot_c:=\{\epsilon_i-\epsilon_j\}$ and the non compact ones 
$\Rgot_n=\{\pm(\epsilon_i+\epsilon_j)\}$. We work with the subsets of positive roots $\Rgot^+_c:=\{\epsilon_i-\epsilon_j, i<j\}$ and $\Rgot^+_n:=\{\epsilon_i+\epsilon_j\}$. The Weyl chamber $\tgot_+\subset \tgot$ 
is defined by the relations  $\langle \alpha,\mu\rangle \geq 0$, $\forall \alpha\in\Rgot^+_c$, namely $\mu_1\geq \cdots\geq \mu_n$. 
The subchamber $\Ccal_n\subset \tgot_+$ is defined by the conditions $\langle \beta,\mu\rangle > 0$, $\forall \beta\in\Rgot^+_n$. Thus 
$X(\mu)\in \Ccal_n$ if and only if $\mu\in\Rbb^n_{++}$.

 If $M\in \spgot(\Rbb^{2n})$, we denote by $\Ocal_M:=\{gM g^{-1}, g\in Sp(\Rbb^{2n})\}$ the corresponding adjoint orbit. 

\begin{lem}
\begin{enumerate}
\item $M\in\mathbf{C}(n)^0$ if and only if there exists $X\in\Ccal_n$ such that $M\in \Ocal_{X}$.
\item Let $\mu\in \Rbb^n_{++}$, and $M\in \Ocal_{X(\mu)}$. The symplectic eigenvalues of the positive definite quadratic form $q(v)= {}^t v J_n M v=\Omega_n(Mv,v)$ are the 
$\mu_1\geq\cdots\geq \mu_n>0$.
\end{enumerate}
\end{lem}

{\em Proof:}  The first point is a classical fact \cite{Vinberg80,Paneitz83}. If $M=gX(\mu)g^{-1}$ with $g\in Sp(\Rbb^{2n})$, we see that 
$$
\Omega_n(Mv,v)=\Omega_n(X(\mu)g^{-1}v,g^{-1}v)=\sum_{k=1}^n \mu_k(v_{k}^2+v_{k+n}^2)
$$
where each $v_j$ is the $j$-th coordinate of the vector $g^{-1}v$. $\Box$

\begin{rem}
In \cite{pep-JIMJ}, we call the interior $\mathbf{C}(n)^0$ of $\mathbf{C}(n)$ the holomorphic cone, since any coadjoint orbit $\Ocal_X\subset \mathbf{C}(n)^0$ admits a canonical structure of 
a K\"{a}hler manifold with a holomorphic action of $K$. These orbits are closely related to the holomorphic discrete series representations of the symplectic group $Sp(\Rbb^{2n})$.
\end{rem}

Thanks to the previous Lemma, we see that the symplectic Horn cone admits the alternative definition:
$$
\SHorn(n)=\left\{(x,y,z)\in(\Rbb_{++}^n)^3\, |\, \Ocal_{X(z)}\subset\Ocal_{X(x)}+\Ocal_{X(y)}\right\}.
$$ 

In the next section, we explain the result of \cite{pep-JIMJ} concerning the determination of $\SHorn(n)$.

\section{Convexity results}\label{sec:convex}

The trace on $\glgot(\Rbb^{2n})$ provides an identification between $\spgot(\Rbb^{2n})$ and its dual $\spgot(\Rbb^{2n})^*$:  to $X\in\spgot(\Rbb^{2n})$ we associate 
$\xi_X\in\spgot(\Rbb^{2n})^*$ defined by $\langle \xi_X,Y\rangle =-\trace(XY)$. Through this identification the causal cone $\mathbf{C}(n)$ becomes 
$$
\mathbf{\widetilde{C}}(n):=\left\{\xi\in\spgot(\Rbb^{2n})^*;  \langle \xi, {\rm Ad}(g)z\rangle \geq 0,\, \forall g\in Sp(\Rbb^{2n})\right\}
$$
where $z=\frac{-1}{2}J_n$. The identification $\spgot(\Rbb^{2n})\simeq\spgot(\Rbb^{2n})^*$ induces several identifications $\kgot\simeq\kgot^*$, $\tgot\simeq\tgot^*$ and $\tgot_{+}\simeq\tgot^*_+$. 
In the latter cases the identifications are done through an invariant scalar product $(-,-)$ on $\kgot^*$.The subchamber $\widetilde{\Ccal}_n\subset \tgot^*_+$ is defined by the conditions: 
$(\alpha,\xi) \geq 0$, $\forall \alpha\in\Rgot^+_c$, and $(\beta,\xi)> 0$, $\forall \beta\in\Rgot^+_n$.

Through $\spgot(\Rbb^{2n})\simeq\spgot(\Rbb^{2n})^*$, the symplectic Horn cone becomes 
$$
\horn_{\rm hol}(Sp(\Rbb^{2n})):=\left\{(\xi_1,\xi_2,\xi_3)\in(\widetilde{\Ccal}_n)^3\, |\, \Ocal_{\xi_3}\subset\Ocal_{\xi_1}+\Ocal_{\xi_2}\right\}.
$$
Here we have kept the notations of \cite{pep-JIMJ}.

We have a Cartan decomposition $\spgot(\Rbb^{2n})=\kgot\oplus\pgot$ with
$$
\pgot:=
\left\{
\begin{pmatrix}A& B\\ B& -A  \end{pmatrix},\ {}^t A= A, \ {}^t B= B
\right\}.
$$
We denote by $\pgot^+$ the vector space $\pgot$ equipped with the complex structure ${\rm ad}(z)$ and the compatible symplectic structure 
$\Omega_{\pgot^+}(Y,Y'):=-{\rm Tr}(J_n[Y,Y'])$: here $\Omega_{\pgot^+}(Y,[z,Y])>0$ for any $Y\neq 0$. 

The action of maximal compact subgroup $K\subset Sp(\Rbb^{2n})$ on 
$(\pgot^+,\Omega_{\pgot^+})$ is Hamiltonian with moment map
$$
\Phi_{\pgot^+}:\pgot^+\to \kgot^*
$$ 
defined by $\langle\Phi_{\pgot^+}(Y),X\rangle=\tfrac{1}{2}\Omega_{\pgot^+}([X,Y],Y)$. If $Y=\begin{pmatrix}A& B\\ B& -A  \end{pmatrix}$, we see that
$\langle\Phi_{\pgot^+}(Y),J_n\rangle={\rm Tr}(A^2+B^2)=\tfrac{1}{2}\|Y\|^2$. Hence the moment map $\Phi_{\pgot^+}$ is proper.


\medskip

We consider the following action of the group $K^3$ on the manifold $K\times K$: 
$$
(k_1,k_2,k_3)\cdot (g,h)=(k_1gk_3^{-1},k_2 hk_3^{-1})
$$

The action of $K^3$ on the cotangent bundle $N:=T^*(K\times K)$ is Hamiltonian with moment map $\Phi_N: N\to \kgot^*\times\kgot^*\times\kgot^*$ 
defined by the relations\footnote{We use the identification $T^*K\simeq K\times \kgot^*$ given by left translations.}
$$
\Phi_N(g_1,\eta_1;g_2,\eta_2)=(-g_1\eta_1, -g_2\eta_2, \eta_1+\eta_2).
$$
Finally we consider the Hamiltonian $K^3$-manifold  $N\times \pgot^+$, where $\pgot^+$ is equipped with the symplectic structure $\Omega_{\pgot^+}$. The action is defined by the relations: $(k_1,k_2,k_3)\cdot (g,h,X)=(k_1gk_3^{-1},k_2 hk_3^{-1},k_3X)$. Let us denote by $\Phi: N\times \pgot^+\to \kgot^*\times\kgot^*\times\kgot^*$ 
the moment map relative to the $K^3$-action~:
\begin{equation}\label{eq:momentcotangent}
\Phi(g_1,\eta_1;g_2,\eta_2,Y)=(-g_1\eta_1, -g_2\eta_2, \eta_1+\eta_2+\Phi_{\pgot^+}(Y)).
\end{equation}

Since $\Phi$ is proper map, the Convexity Theorem \cite{Kirwan84,LMTW} tell us that 
$$
\Delta(N\times \pgot^+):={\rm Image}(\Phi)\bigcap \, \tgot_{+}^*\times \tgot_{+}^*\times \tgot_{+}^*
$$
is a closed, convex, and locally polyhedral set.

The map $\mu\mapsto X(\mu)$ defines an isomorphism of $\Rbb^n$ with $\tgot\simeq\tgot^*$ that induces an identification of $\Rbb^n_{++}$ with $\Ccal_n\simeq \widetilde{\Ccal}_n$.
Recall that on $\tgot^*\simeq\Rbb^n$, we have a natural involution that sends $\mu=(\mu_1,\ldots,\mu_n)$ to $\mu^*:=(-\mu_n,\ldots,-\mu_1)$. The following result is proved in  \cite{pep-JIMJ} (see Theorem B). 

\medskip

\begin{theo}\label{theo:pep}
An element $(x,y,z)\in(\Rbb_{++}^n)^3$ belongs to $\horn_{\rm hol}(Sp(\Rbb^{2n}))$ if and only if 
$$
(x,y,z^*)\in \Delta(N\times \pgot^+).
$$
\end{theo}

\medskip

Recall that a Hermitian matrix $M$ majorizes another Hermitian matrix $M$ if $M-M'$ is positive semidefinite (its eigenvalues are all nonnegative). In this case, we write $M\geq M'$. 

\begin{prop}\label{prop:delta}
Let $(x,y,z)\in(\Rbb_{+}^n)^3$. Then  $(x,y,z^*)\in \Delta(N\times \pgot^+)$ if and only if there exist Hermitian matrices $A,B,C$ such that \ 
$\mathrm{s}(A)=x$, $\mathrm{s}(B)=y$, $\mathrm{s}(C)=z$ and $C\geq A+B$.
\end{prop}

{\em Proof:}  The map $\begin{pmatrix}A& -B\\ B& A  \end{pmatrix}\mapsto A-iB$ defines an isomorphism between $K$ and the unitary group $U(n)$. 
Let us denoted by $S^{2}(\Cbb^n)$  the vector space of complex $n\times n$ symmetric matrices that is equipped with 
the following action of $U(n)$: $k\cdot M= kM{}^t\!k$. The map $\begin{pmatrix}A& B\\ B& -A  \end{pmatrix}\mapsto A-iB$ defines an 
isomorphism between the $K$-module $\pgot^+$ and the $U(n)$-module 
$S^{2}(\Cbb^n)$. Through this identifications the moment map $\Phi_{\pgot^+}:\pgot^+\to \kgot^*$ 
becomes the map $\Phi_{S^2}:S^{2}(\Cbb^n)\to \ugot(n)$ defined by the relations
$$
\Phi_{S^2}(M)=-2iM\overline{M}.
$$

So we know that the moment polytope $\Delta$ relative to the Hamiltonian action of $U(n)^3$ on $T^* U(n)\times T^* U(n) \times S^{2}(\Cbb^n)$ is equal to 
$\Delta(N\times \pgot^+)$. A small computation shows that $(x,y,z^*)\in \Delta$ if and only if there exists Hermitian matrices $A,B,C$ and $M\in S^{2}(\Cbb^n)$ such that 
$$
\mathrm{s}(A)=x,\hspace{3mm} \mathrm{s}(B)=y, \hspace{3mm}\mathrm{s}(C)=z\quad \text{and}\quad A+B+2M\overline{M}=C.
$$
The existence of $M\in S^{2}(\Cbb^n)$ satisfying the condition $A+B+2M\overline{M}=C$ is equivalent to $C\geq A+B$. The proof is then completed. $\Box$

\medskip

S. Friedland \cite{Friedland} considered the following question: {\em which eigenvalues $(\mathrm{s}(A),\mathrm{s}(B),\mathrm{s}(C))$ can occur if $C\geq A+B$}. His solution was in terms of linear inequalities, 
which includes Klyachko’s inequalities, a trace inequality and some additional inequalities. Later, W. Fulton \cite{Fulton-00-bis} proved the additional inequalities are unnecessary. Let us summarize their result in the following Theorem.

\begin{theo}[\cite{Friedland,Fulton-00-bis}]\label{theo:FF}A triple $x,y,z\in\Rbb^n_+$ occurs as the eigenvalues of $n$ by $n$ 
Hermitian matrices $A$, $B$, $C$ with $C\geq A+B$ if and only it satisfies $|x|+|y|\leq |z|$ and $(\star)_{I,J,K}$ for all $(I,J,K)$ of cardinality $r<n$ such that  $c_{IJ}^{K}=1$.
\end{theo}

\medskip

The combination of Theorems \ref{theo:pep} and \ref{theo:FF} with Proposition \ref{prop:delta} completes the proof of Theorem \ref{theo:main}.

\bibliographystyle{crplain}

\end{document}